\documentclass[12pt]{article}

\usepackage{amsfonts}

\usepackage{amsmath}
\usepackage{graphicx}

\newtheorem{theorem}{Theorem}

\newtheorem{proposition}{Proposition}

\newcommand{\pibar}{\overline{\Pi}}

\newcommand{\veps}{\varepsilon}

\newcommand{\eqdr}{\stackrel{\mathrm{D}}{=}}
\newcommand{\R}{\mathbb{R}}
\newcommand{\halmos}{$\sqcup\!\!\!\!\sqcap$}

\newcommand{\be}{\begin{equation}}
\newcommand{\ee}{\end{equation}}
\newcommand{\ben}{\begin{equation*}}
\newcommand{\een}{\end{equation*}}
\newcommand{\ba}{\begin{aligned}}
\newcommand{\ea}{\end{aligned}}
\newcommand{\rmd}{{\rm d}}
\newcommand{\rmi}{{\rm i}}
\newcommand{\wh}{\widehat}

\newcommand{\Xbar}{\overline{X}}

\newcommand{\Gtau}{G_{\tau_u-}}
\newcommand{\pard}{\frac{\partial_-}{\partial_- u}}

\newcommand{\KeesT}{\rule{0pt}{5.6ex}}
\newcommand{\KeesB}{\rule[-5ex]{0pt}{0pt}}
\newenvironment{narrow}[2]{%
  \begin{list}{}{%
    \setlength{\topsep}{0pt}%
    \setlength{\leftmargin}{#1}%
    \setlength{\rightmargin}{#2}%
    \setlength{\listparindent}{\parindent}%
    \setlength{\itemindent}{\parindent}%
    \setlength{\parsep}{\parskip}%
  }%
  \item[]
}{\end{list}}

\numberwithin{equation}{section}
\numberwithin{lemma}{section}
\numberwithin{proposition}{section}
\numberwithin{theorem}{section}
\numberwithin{corollary}{section}
\begin{document}

\date{}
\title{Asymptotic Distributions of the
Overshoot and Undershoots for the L\'{e}vy Insurance Risk Process
in the Cram\'er and  Convolution Equivalent Cases}
\author{Philip S. Griffin, Ross A. Maller and Kees van Schaik\thanks{Research
partially supported by ARC Grant DP1092502}
\\Syracuse University, the Australian National University
\\ and the University of Manchester}
\maketitle

\begin{abstract}
\noindent
Recent models of the insurance risk process use a L\'evy process to generalise the traditional
Cram\'er-Lundberg compound Poisson model.
This paper is concerned with the behaviour of the distributions of the overshoot and undershoots of a high level,
for a L\'{e}vy process which drifts to $-\infty$ and satisfies a Cram\'er or a convolution equivalent condition.
We derive these asymptotics under minimal conditions in the
Cram\'er case, and compare them with known results for the convolution equivalent case,
drawing attention to the striking and unexpected fact that they become identical when certain
parameters tend to equality.
Thus, at least regarding these quantities, the  ``medium-heavy" tailed  convolution equivalent model segues into
the ``light-tailed" Cram\'er model in a natural way.
This suggests a usefully expanded flexibility for modelling the insurance risk process.
We illustrate this relationship by comparing the asymptotic distributions obtained for the overshoot
and undershoots, assuming the L\'evy process belongs to the ``GTSC" class.
\end{abstract}

\noindent\textit{Corresponding Author:}
\\ Ross A. Maller
\\ Mathematical Sciences Institute
\\ The Australian National University
\\ PO Canberra ACT
\\ Australia
\\ Ross.Maller@anu.edu.au



\bigskip
\noindent\textit{Insurance Branch Category Codes:} IM11, IM13.

\vspace*{10pt} \setcounter{equation}{0}
\section{Introduction}\label{s1}

There has recently been a great deal of interest in the insurance and related literature
concerning continuous time risk models, analysed under  assumptions which have been found to be realistic
for classical discrete time (compound Poisson) models.
A very tractable continuous time generalisation of the classical model is to use a L\'evy process,
$(X_t)_{t\ge 0}$, for the risk process.
Results for interesting quantities in this model (overshoot and undershoots)
can be  obtained under an analogue of the classical Cram\'er-Lundberg condition,
as we do in this paper
(see also Biffis and Kyprianou (2010), Schmidli (1995), and their references for related discussions),
or under a  ``convolution equivalent" condition
on the tail of the L\'evy  measure of $X$ (see Kl\"{u}ppelberg,  Kyprianou and Maller (2004),
Tang and Wei (2010), Embrechts and Veraverbeke (1982), and their references).
Our aim in the present paper is to review the L\'evy setup,
derive some new results for the Cram\'er case,
and, further, to draw out the remarkable connection between the Cram\'er and convolution equivalent results
 mentioned in the abstract.
This provides a very useful and practical synthesis of the two different approaches.

The L\'evy setup is as follows.
Suppose that $X=\{X_{t}: t \geq 0 \}$,
$X_0=0$,  is a real valued L\'{e}vy process
defined on $(\Omega, {\cal F}, P)$,
with triplet $(\gamma, \sigma^2, \Pi_X)$,
$\Pi_{X}$ being the L\'{e}vy measure of $X$, $\gamma\in \R$, $\sigma^2\ge 0$.
Thus the characteristic function of $X$
is given by the L\'{e}vy-Khintchine
representation, $E(e^{i\theta X_{t}}) = e^{t \Psi(\theta)}$,
where
\begin{equation}\label{lrep}
\Psi(\theta) =
 \rmi\theta \gamma - \sigma^2\theta^2/2+
\int_{\R\setminus\{0\}}(e^{\rmi\theta x}-1-
\rmi\theta x \mathbf{1}_{\{|x|<1\}})\Pi_X(\rmd x),
\ {\rm for}\  \theta \in \mathbb{R}.
\end{equation}

We will consider limiting distributions of
the overshoot/undershoot of the process above/below a high level,
and the undershoot from the previous maximum,
viz, the quantities $X_{\tau_u}-u$, $u-X_{\tau_u-}$,
and $u-\overline X_{\tau_u-}$,
conditional on $\tau_u < \infty$, where
$\overline{X}_{t} = \sup_{0\le s\leq t}X_{s}$, and
$\tau_u$
is the first passage time above level $u>0$:
\begin{eqnarray*}
\tau_u= \inf \{t \geq0 : X_{t}>u \},~~ u > 0
\end{eqnarray*}
(with $\tau_u=\infty$ if $X_t\le u$ for all $t>0$).
We assume either the  Cram\'er condition, namely, that
$Ee^{\nu_0X_1}=1$ for some $\nu_0>0$,
or, alternatively, $Ee^{\alpha X_1}<1$ for some $\alpha>0$  together with a  ``convolution equivalent" condition (see Section \ref{s5})
on the tail of the L\'evy   measure of $X$.
Results under the latter condition were derived in
Kl\"{u}ppelberg et al.  (2004),
Doney and Kyprianou (2006), Park and Maller (2008),
and Griffin and Maller (2010),
which papers we will draw on for the relevant comparisons.
The asymptotic behaviour of the process is quite different in these two settings, but we will make explicit some remarkable parallels
between the associated overshoot and undershoots.

We assume that the process $X$ drifts to $-\infty$, modelling the situation of
an insurance risk process with premiums and other income producing a
downward drift in $X$, and claims represented by positive
jumps.  Thus the process
$X$, called the {\it claim surplus process}, represents the excess in
claims over income.
We think of an insurance company starting with a positive reserve
$u$, with  ruin  occuring
if this level is exceeded by $X$.
We will refer to this as the {\it General L\'evy Insurance Risk Model}.
It is a generalisation of the classical {\it Cram\'{e}r-Lundberg model},
which arises in the special case
when the claim surplus process is taken to be
\be\label{CL}
X_t=\sum_{i=1}^{N_t} U_i - rt, \ r>0,
\ee
where the $(U_i)_{i=1,2,\ldots}$ are i.i.d. positive random variables  independent
of the  Poisson process $(N_t)_{t\ge 0}$.
Here $r$ represents the rate of premium
inflow and $U_i$ the size of the $i$th claim.

The general model allows for income other than through premium inflow and a more realistic claims structure.
Large claims, modeled by large jumps, are offset by premiums collected at a roughly constant rate.
These claims, spaced out by independent exponential times, correspond
reasonably well to ``disasters", of which there have been many in recent history.
The compensated small jumps, which may be countably infinite in number,
correspond to ongoing minor claims, with
their compensator understood as the aggregate of premiums required to
offset the high intensity of claims.
The assumption $\lim_{t\to\infty} X_t= -\infty$ a.s. is a reflection of premiums being
set to avoid almost sure ruin.

We say that ``ruin occurs" if  $\tau_u<\infty$.
Interest centres on the properties of the process when ruin occurs,
as a kind of worse case scenario, so
results in this area are stated as limit theorems
conditional on $\tau_u<\infty$.
The assumption $X_t\to -\infty$ a.s. implies
$P(\tau_u<\infty)<1$,  so ruin is not certain in this model,
while  $P(\tau_u<\infty)\to 0$
as the initial level $u\to\infty$,
so ruin becomes less and less probable as the initial capital
is increased, as is logical.
But $P(\tau_u<\infty)>0$ for finite $u$,
and it is convenient to define by elementary means
a new probability measure $P^{(u)}$ given by
\begin{equation}\label{pudef}
P^{(u)}(\ \cdot\ )=P(\ \cdot\ |\tau_u<\infty).
\end{equation}
Our convergence results will relate to this distribution.
In the Cram\'er case they rely on some fundamental methods introduced by
Bertoin and Doney (1994), and, like them, we will not need further
assumptions on the tail behaviour, or otherwise, of $X$.

In what follows, Section \ref{s2} lays out the basic foundational assumptions for the model, and Section \ref{s3} lists the marginal distributions we will need, calculated from a ``quintuple law".
We separate the setup and results
for the Cram\'er  case in Section \ref{s4}
and for the convolution equivalent case  in Section \ref{s5}.
Section \ref{s6} gives a summary of the general Cram\'er and convolution equivalent results, and in Section \ref{s7} we give a diagrammatic representation and comparison of the  results for a particular class (the ``GTSC" class) of L\'evy processes.
A discussion of the significance of the results is in Section \ref{s8}.
Proofs are collected in Section \ref{s9}.

\setcounter{equation}{0}
\section{Fluctuation Foundations}\label{s2}

A L\'evy process has certain ``fluctuation" quantities associated with it as follows.
Let $L_t$ denote the local time of $X$ at its maximum  and
$(L^{-1}_t,H_t)_{t \geq 0}$
the ascending bivariate ladder process of $X$,
which is defective under {\em the condition}
$\lim_{t\to \infty} X_{t} = -\infty$ a.s. {\em (which will always
obtain in this paper)}.
The nondefective descending ladder process
is denoted by $(\wh L^{-1}_t, \widehat{H}_{t})_{t\geq0}$ (taking
{\em positive} values).

The defective process
$(L^{-1},H)$ is obtained from  a nondefective process
$({\cal L}^{-1}, {\cal H})$ by exponential killing.
To be precise, $L_\infty$ has an exponential distribution with some parameter $q > 0$ say,
and
$\left((L^{-1},H):t<L_\infty\right)
\eqdr\left(({\cal L}^{-1},{\cal H}):t<e(q)\right)$,
where $e(q)$ is independent of $({\cal L}^{-1},{\cal H})$ and has
exponential distribution with parameter $q$.
Processes
$({\cal L}^{-1}, {\cal H})$ and $(\wh L^{-1}, \widehat{H})$
are nondefective bivariate subordinators, and the marginal processes are
nondefective univariate subordinators.

Denote the bivariate L\'{e}vy measure of
$({\cal L}^{-1},{\cal H})$
by $\Pi_{{\cal L}^{-1},{\cal H}}(\cdot,\cdot)$.
The Laplace exponent $\kappa(a,b)$ of  $(L^{-1},H)$,
which satisfies
\begin{equation} \label{kapdef}
e^{-\kappa(a,b)}=
E (e^{-a{L}^{-1}_1 -b{H}_1};L_\infty >1)=
E( e^{-a{\cal L}^{-1}_1 -b{\cal H}_1}; e(q)>1)
\end{equation}
for values of $a,b\in \R$ for which the expectation is finite,
may be written
\begin{eqnarray} \label{kapexp}
\kappa(a,b)
&=&
q+\rmd_{L^{-1}}a+\rmd_Hb+\int_{t\ge0}\int_{x\ge0}
\left(1-e^{-at-bx}\right)
\Pi_{{\cal L}^{-1},{\cal H}}(\rmd t, \rmd x),
\end{eqnarray}
where $\rmd_{L^{-1}}\ge 0$ and $\rmd_H\ge 0$ are drift constants.
Denote the marginal L\'{e}vy measures of ${\cal L}^{-1}$ and ${\cal H}$
(equivalently of $L^{-1}$ and $H$)
by $\Pi_{{\cal L}^{-1}}$ and $\Pi_{\cal H}$, and similarly for $\wh L^{-1}$ and $\wh H$.

The bivariate renewal function of $(L^{-1},H)$ is
 \begin{equation}\label{Vkdef}
V(t,x)= \int_0^\infty P(L_s^{-1}\le t,H_s\le x)\rmd s
= \int_0^\infty e^{-qs} P({\cal L}^{-1}_s\le t,{\cal H}_s\le x)
\rmd s, \ {\rm for}\ t\ge 0, x \ge 0.
\end{equation}
Its  Laplace transform is given by
\begin{equation}\label{Vkap}
\int_{t\ge 0}\int_{x\ge 0}e^{-at-bx} V(\rmd t,\rmd x)=\int_{s\ge 0}
E(e^{-a L_s^{-1}-bH_s}; L_\infty >s)\rmd s =\frac{1}{\kappa(a,b)}
 \end{equation}
for values $a,b$ for which $\kappa(a,b)>0$.
The marginal renewal function for $H$ is
\begin{equation}\label{VHdef}
V_H(x)= \int_0^\infty P(H_s\le x)\rmd s
= \int_0^\infty e^{-qs} P({\cal H}_s\le x)\rmd s
=\lim_{t\to\infty}V(t,x).
\end{equation}
Analogously the renewal function of $\wh H$ will be denoted $\wh V_{\wh H}$.
The  tails of $\Pi_X$ are
\begin{eqnarray*}
\overline{\Pi}_X^{+}(x) = \Pi_X\{(x,\infty)\}, \
\overline{\Pi}_X^-(x)= \Pi_X\{(-\infty, -x)\},
\ {\rm and}\
\overline{\Pi}_X(x) =\overline{\Pi}_X^+(x) +
\overline{\Pi}_X^-(x), \ x > 0,
\end{eqnarray*}
with similar notations $ \overline{\Pi}_{\cal H}(x)$ and $
\overline{\Pi}_{\wh H}(x)$ for ${\cal H}$ and $\wh H$.

\setcounter{equation}{0}
\section{Distributional Identities}\label{s3}
In this section we present some general formulae for distributions
of the ruin variables (overshoot and undershoots), in terms
of the fluctuation variables.
From these we can work out the asymptotic distributions of the ruin variables.

The formulae are based on the ``quintuple law" of Doney and Kyprianou (2009, Theorem 3),
which gives the joint distribution of three position  and two time variables associated with ruin.
Rather than state the quintuple law in its general form, which will not be used in this paper,
we give a formula for the marginal trivariate distribution of the overshoot and undershoots which
can be calculated directly from it; namely, we have,
for $u>0$, $x\ge 0$, $v\ge 0$, $0\le y\leq u \wedge v$,
\begin{eqnarray}\label{cdf4}
&&P\left(X_{\tau_u}-u \le x, u-X_{\tau_u-} \le v,
u-\overline{X}_{\tau_u-} \le y,\tau_u<\infty\right)
\nonumber\\
&=&
\int_{0\le y'\le y}
\int_{0\le v'\le v-y'}
\Pi_X\{(v'+y',v'+y'+x]\}
\widehat{V}_{\wh H}(\rmd v')|V_H(u-\rmd y')|
+\rmd _H V'_H(u).
 \nonumber\\
\end{eqnarray}

Some L\'evy processes possess the property of ``creeping"
over level $u$,
meaning that $X_{\tau_u}=u$ can occur with positive probability.
(This cannot happen in the classical Cram\'er-Lundberg model.)
By a well known result of Kesten, see Bertoin (1996, Theorem VI.19),
\be\label{kes}
P(X_{\tau_u}=u)=\rmd _H V'_H(u),
\ee
with the understanding that $\rmd _H V'_H(u)=0$ when $\rmd _H=0$
($V_H$ need not be differentiable  in this  case).
The second term in \eqref{cdf4} is a consequence of the possibility of $X$ creeping
over level $u$.

With the help of Vigon's (2002) ``equation amicale invers\'ee",
viz
\begin{equation}\label{DKcor6}
\pibar_{\cal {H}}(x)
=\int_{v\ge 0}\pibar_X^+(x+v) \wh V_{\wh H}(\rmd v),\quad  x> 0,
\end{equation}
the marginal distributions of the position variables can
be written down as follows.

\bigskip \noindent  {\bf Marginal Distribution of the Overshoot:} for $u>0$, $x\ge 0$,
\begin{eqnarray}\label{om}
P\left(X_{\tau_u}-u \le x, \tau_u<\infty\right)
=
\int_{[0,u]}
\left(\pibar_{{\cal H}}(u-y)-\pibar_{{\cal H}}(u+x-y)\right)V_H(\rmd y)
+\rmd _H V'_H(u).
\nonumber\\
\end{eqnarray}
Letting $x\to\infty$ in this gives
 \begin{equation} \label{tm}
P\left(\tau_u<\infty\right)
=\int_{[0,u]}\pibar_{{\cal H}}(u-y)V_H(\rmd y)+\rmd _H V'_H(u).
 \end{equation}

\bigskip \noindent  {\bf Marginal Distribution of the Undershoot:}
Define
\begin{equation} \label{gdef}
g_v(y):=\int_{z \in[0, v-y]}
\overline{\Pi}_{X}^{+}(z+y)\wh V_H(\rmd z),\ 0\le y\le v.
\end{equation}
Then for $u>0$, $v\ge 0$,
\begin{eqnarray} \label{um}
P(u-X_{\tau_u-} \le v,\tau_u<\infty)
=\int_{u-(u\wedge v)\le y\le u}g_v(u-y) V_H(\rmd y)
+\rmd _H V'_H(u).
\end{eqnarray}
For later reference we note that $g_v(y)$ is decreasing $y$, increasing in $v$, and by \eqref{DKcor6} we have
\be\label{ginf}
g_v(y)\uparrow\pibar_{{\cal H}}(y)\ \text{ as $v\uparrow \infty$,   all $y>0$}.
\ee
Recall that for any subordinator $\cal{H}$,
\be\label{sub0}
y\pibar_{{\cal H}}(y)\to 0\ \text{ as $y\to 0$,}
\ee
thus while it is possible that $g_v(0)=\infty$, by \eqref{ginf} and \eqref{sub0} we have
\be\label{ygto0}
y g_v(y)\to 0 \quad \text{as } y\to 0.
\ee

\bigskip \noindent  {\bf Marginal Distribution of the Undershoot
from the Previous Maximum:}
 for $u>0$, $0\le y\le u$,
\begin{eqnarray} \label{3.7}
P\left(u-\overline X_{\tau_u-} \le y,\tau_u<\infty \right)
=\int_{u-y\le z\le u}\overline{\Pi}_{{\cal H}}(u-z)V_H(\rmd z)
 +\rmd _H V'_H(u).
\end{eqnarray}
Of course this undershoot lies between 0 and $u$.
\medskip

By incorporating an extension of the quintuple law which can be
found in Griffin and Maller (2011),
similar, but slightly more complicated formulae, can be written down for the distributions of the
ruin time and the time of the last maximum of $X$ (i.e., the time of minimum surplus) prior to ruin.
These will not be used in this paper, but we record them here for completeness.  Let
${G}_{t} = \sup \{s\leq t : X_{s}=\overline{X}_{s}\}$
be the time of the last maximum of $X$ up to time $t$.
Thus $\Gtau$ is the time of the last maximum of $X$ prior to ruin.
We then have the following analogue of \eqref{cdf4} for the time variables:
 for $u>0$, $s\ge 0$, $t\ge 0$,
\begin{eqnarray}\label{cdf6a}
&&P\left(\tau_u-\Gtau \le s, \Gtau \le t\right)
=
\int_{0\le y\le u}
\Pi_{{\cal L}^{-1},{\cal H}}([0,s],(y,\infty))|V(t,u-\rmd y)|
 \nonumber\\
&&
\qquad \qquad \qquad \qquad
\qquad \qquad \qquad \qquad \qquad
+\rmd_H \pard V(t,u),
\end{eqnarray}
where $\partial_-/\partial_-u$ denotes left derivative.
The second term in \eqref{cdf6a} is a consequence of the possibility of creeping by time $t$;
see Griffin and Maller (2011) for details.

From this we obtain the

\bigskip \noindent  {\bf Marginal Distribution of the Time of Minimum Surplus Prior to Ruin:}
for $u>0$, $t\ge 0$,
\begin{eqnarray}\label{cdf6}
&&P\left(\Gtau \le t, \tau_u<\infty\right)
=
\int_{0\le y\le u}
\pibar_{{\cal H}}(y)|V(t,u-\rmd y)|
+\rmd_H \pard V(t,u).
\end{eqnarray}

\bigskip \noindent  {\bf Marginal Distribution of the Ruin Time:}
for $u>0$, $s\ge 0$,
\begin{equation}\label{taudis}
P\left(\tau_u\le s\right)
=\int_{0\le y\le u}\int_{0\le t \le s}
\Pi_{{\cal L}^{-1},{\cal H}}([0,s-t],(y,\infty)) |V(\rmd t,u-\rmd y)|
+\rmd _H \pard V(s,u).
\end{equation}

\setcounter{equation}{0}
\section{The  Cram\'er Setup and Results}\label{s4}

Our assumption throughout  this section
is the {\em Cram\'er condition}, namely, that
\begin{equation}\label{cramer}
Ee^{\nu_0X_1}=1, \ {\rm  for\ some}\  \nu_0>0.
\end{equation}
That this implies
the process drifts to $-\infty$, i.e.,
$\lim_{t\to\infty}X_{t} = -\infty$ a.s., and some further
properties are summarised in the next proposition. Its proof is
well known and omitted.

\begin{proposition}\label{t1}
Assume  \eqref{cramer}. Then $EX_1$ is well defined,
and $EX_1\in [-\infty,0)$, and  so
$\lim_{t \to \infty} X_t=-\infty$ a.s.
Further, $ Ee^{\nu X_1}$ is finite and nonzero for all
$\nu\in[0,\nu_0]$, and $ Ee^{\nu X_1}$ increases in $\nu$ for
$\nu\ge \nu_0$; in particular,  $1\le Ee^{\nu X_1}\le \infty$ for
$\nu\ge \nu_0$.
Further, $E(X_1e^{\nu X_1})\in (0,\infty)$
for all $\nu$ in a left neighbourhood $(\nu^*,\nu_0)$, $\nu^*<\nu_0$,  of $\nu_0$,
and  $E(X_1e^{\nu_0 X_1})>0$ (possibly, $E(X_1e^{\nu_0 X_1})=+\infty$).
\end{proposition}

Additionally, use of the Wiener-Hopf factorisation as in the proof of
Prop. 5.1 in Kl\"uppelberg et al. (2004), shows that
 \begin{equation}\label{did}
E(e^{\nu_0 X_1})=1\  {\rm implies}  \ E(e^{\nu_0H_1};L_\infty >1)=1.
\end{equation}
Consequently,  under \eqref{cramer},
it follows from \eqref{kapexp} with $a=0$ and $b=-\nu_0$ that
\begin{equation}\label{qid}
 q=\nu_0\rmd_H+\int_{h>0}\left(e^{\nu_0 h}-1\right) \Pi_{{\cal H}}(\rmd h).
\end{equation}

Bertoin and Doney (1994) show that Cram\'er's estimate for ruin continues to hold for general L\'evy processes
under \eqref{cramer}.
To state their result recall that $X$ is a compound Poisson process if $\Pi_X$ is finite and $X$ is the sum of its jumps.
Their main result (see also Theorem 7.6 of Kyprianou (2005) and Section XIII.5 of Asmussen (2002)) may then be stated as follows:
{\it Suppose the support of $\Pi_X$ is non-lattice in the case that
$X$ is compound Poisson. Then
 \begin{equation}\label{BDmain}
\lim_{u\to\infty}e^{\nu_0 u}P(\tau_u<\infty)=\frac{q}{\nu_0m^*},
 \end{equation}
where}
 \begin{eqnarray}\label{rid}
m^*:=
\rmd_H+\int_{h>0}h e^{\nu_0 h}\Pi_{\cal H}(\rmd h)\in (0,\infty].
 \end{eqnarray}
A similar result holds in the lattice case if the limit is taken through points in the lattice span,
but to avoid repetition we will assume hence forth the non-lattice version  of \eqref{BDmain}.

A further useful result in Bertoin and Doney (1994) is the following version of the renewal theorem for L\'evy processes. Let

\begin{equation}\label{tauV}
U^*(x)=\int_{y\le x}e^{\nu_0 y}V_H(\rmd y).
 \end{equation}
Then Bertoin and Doney (1994) show that $U^*$ is a renewal function satisfying
\begin{equation}\label{U*con}
\lim_{u\to\infty}
\left(U^*(u+x)-U^*(u)\right)=\frac{x}{m^*}, \ x\in \R,
\end{equation}
where the righthand side of \eqref{U*con}
is interpreted as zero if $m^*=\infty$.  Using \eqref{qid} and \eqref{Vkap} with $a=0$, it is easily seen that in fact $U^*$ is the renewal function of the nondefective subordinator $H^*$ which has drift $d_{H^*}$ and L\'evy measure $\Pi_{H^*}$ given by
\be\label{dpstar}
d_{H^*}=d_{H}, \quad\text{and}\quad  \Pi_{H^*}(\rmd y)=e^{\nu_0 y}\Pi_{H}(\rmd y).
\ee

In light of \eqref{BDmain} and \eqref{U*con}, it is important to know when $m^*<\infty$. From Theorem 8 of  Doney and Maller (2002b), for example, it follows that
\begin{equation}\label{mu*1}
m^*<\infty \text{ if and only if } EX_1e^{\nu_0X_1}<\infty.
 \end{equation}
Thus in the Cram\'er set up we will also assume
 \begin{equation}\label{mu*}
EX_1e^{\nu_0X_1}<\infty.
 \end{equation}
 We note from \eqref{rid} that one immediate consequence of \eqref{mu*} is
\begin{equation}\label{picon}
 \lim_{u\to\infty} ue^{\nu_0 u} \pibar_{\cal H}(u)=0.
  \end{equation}

The following theorem is  proved in Section \ref{s9}.

\begin{theorem}\label{th1}
Assume \eqref{cramer} and \eqref{mu*}.
Then we have, for $x\ge 0$,
\begin{equation}\label{t2a}
\lim_{u \to \infty}
P^{(u)}\left(X_{\tau_u}-u\le x\right)=
1-\frac{1}{q}
\int_{(0,\infty)}\left(e^{\nu_0y}-1\right)\Pi_{\cal H}(x+\rmd y),
\end{equation}
 \begin{equation} \label{t5a}
\lim_{u\to\infty}
P^{(u)}\left(u-X_{\tau_u-}\le x\right)=
\frac{\nu_0\rmd_H}{q}
+\frac{1}{q} \int_0^x g_x(y)\rmd\left(e^{\nu_0y}\right),
 \end{equation}
and
 \begin{equation} \label{t6a}
\lim_{u\to\infty}
P^{(u)}\left(u-\overline X_{\tau_u-}\le x\right)=
\frac{\nu_0\rmd_H}{q}
+ \frac{1}{q}\int_0^x
\pibar_{\cal H}(y)\rmd(e^{\nu_0 y}).
 \end{equation}
In particular, setting $x=0$,
\begin{equation}\label{t3a}
\lim_{u \to \infty}
 P^{(u)}\left(X_{\tau_u}=u\right)=
\lim_{u\to\infty}P^{(u)}\left(X_{\tau_u-}=u\right)=
\lim_{u\to\infty}P^{(u)}\left(\overline X_{\tau_u-}=u\right)
=\frac{\nu_0 \rmd_H}{q}.
\end{equation}
\end{theorem}
\medskip

\setcounter{equation}{0}
\section{The  Convolution Equivalent Setup and Results}\label{s5}

The subexponential, or, more
generally, the  convolution equivalent distributions, are used in
modeling heavy tailed data such as often occur in insurance
applications;
we refer to Embrechts and Goldie (1982),  Embrechts et al. (1979), and their references, for discussion and properties.
An important class of distributions
which are convolution
equivalent  for certain values of the parameters
are the generalized inverse Gaussian distributions;  see
Kl\"uppelberg (1989).

We briefly recap the main ideas.
For background and rationale on the convolution equivalent assumptions see Kl\"uppelberg et al. (2004) and Griffin and Maller (2010).
    We say that a distribution $F$
    on $ [0, \infty)$ with tail $\overline{F}$ belongs to
    the {\em class} ${\cal L}^{(\alpha)}$, $\alpha \geq 0$, if
\begin{eqnarray}\label{Lalph}
    \lim_{u \to \infty}
\frac{\overline{F}(u-x)}{\overline{F}(u)} =
    e^{\alpha x}, ~~\mbox{for}~~x \in (-\infty, \infty).
\end{eqnarray}
Tail functions $\overline{F}$ such that $\overline{F}(\log u)$
is regularly varying with index $-\alpha$, $\alpha \geq 0$,
as $u\to \infty$, are in  ${\cal L}^{(\alpha)}$.
With $*$ denoting convolution, a distribution
$F$ is said to be {\em convolution equivalent}, i.e., in the {\em class} ${\cal S}^{(\alpha)}$, $\alpha \geq 0$,
if
$F \in {\cal L}^{(\alpha)}$, and, in addition,
\begin{eqnarray}\label{Salph}
    \lim_{x \to \infty}
\frac{\overline{F^{2*}}(x)}{\overline{F}(x)} :=
2\delta_{\alpha}(F) < \infty,
\end{eqnarray}
where
\[
\delta_{\alpha}(F):= \int_{(0, \infty)} e^{\alpha x}F(\rmd x)
\]
 is the moment generating function of $F$.
The parameter $\alpha$ is referred to as the {\em index} of the
class. When $\alpha=0$, ${\cal S} := {\cal S}^{(0)}$ is the
class of {\em subexponential distributions}.

{\em From now on we keep $\alpha>0$.}

Distributions in the class  ${\cal S}^{(\alpha)}$, for
$\alpha >0$, satisfy $\delta_{\alpha}(F) < \infty$ and
$\delta_{\alpha+\veps}(F) = \infty$ for every $\veps>0$.  They are often referred to as
``near to exponential"  in
the sense that their tails are ``close" to the tail of the exponential distribution with parameter $\alpha$.
For example if
$\overline{F}(x)\sim x^{-1-\rho}e^{-\alpha x}$ for some $\rho>0$ and $\alpha>0$,
then $F \in {\cal S}^{(\alpha)}$; see Chover, Ney and Wainger (1973).
However, note that the exponential distribution itself
is not in any ${\cal S}^{(\alpha)}$.

We can take the tail of any L\'evy  measure, assumed
nonzero on some interval contained in $(0,\infty)$, to be the
tail
of a  distribution function on $[0,\infty)$, after
renormalisation.
With this convention, we say then that the measure is in
$\mathcal{L}^{\left( \alpha \right)}$ or $\mathcal{S}^{\left(
\alpha \right)}$ if this is true for the measure with the
corresponding  (renormalised) tail. The results of
Kl\"{u}ppelberg et al. (2004) are restricted  to the case where the distribution or
measure is not concentrated on a lattice  in $\R$, which we will assume here also,
with the understanding that the alternative can be handled by obvious
modifications.  Beyond this, the main assumptions are that
 for some $\alpha > 0$,
\begin{equation}\label{c1}
\overline \Pi_{X}^{+}(0) > 0, \ \
\overline{\Pi}_{X}^+ \in {\cal S}^{(\alpha)}\ \
\mbox{and}\ \
E e^{\alpha X_1} < 1.
\end{equation}
The case when $\Pi_X$
is concentrated on $(-\infty,0)$, i.e., $X$ is spectrally negative,
is rather easy to handle, so it is  excluded throughout.
(See Remark 4.6, p. 1780, of \cite{kkm}.)

The assumption \eqref{c1}
means that we
analyse the non-Cram\'{e}r case.
This is because if $Ee^{\alpha X_1}<1$  and $Ee^{\nu_0 X_1}=1$, then $0<\alpha<\nu_0$
from  Proposition \ref{t1}.  But membership
in the class  ${\cal S}^{(\alpha)}$
means that $Ee^{(\alpha+\veps) X_1} = \infty$ for every $\veps>0$,
contradicting $Ee^{\nu_0X_1}=1$.
Thus the Cram\'er case is
mutually exclusive to the convolution equivalent case as we consider it here.
The  condition $E e^{\alpha X_1} < 1$ in \eqref{c1} also ensures that $ \lim_{t \to \infty} X_{t}=-\infty$ a.s.

It's convenient to define constants $\beta_1$, $\beta_2$, by
  \begin{equation}\label{betadef}
\beta_1 = -\ln Ee^{\alpha X_1}= -\Psi(-\rmi \alpha)
\ {\rm and}\
\beta_2=\frac{\kappa(0,-\alpha)}{q}
 \end{equation}
 where recall $\kappa$ is given by \eqref{kapexp}.
 Under \eqref{c1},  $\beta_1>0$ and $0<\beta_2<1$.  The first is obvious and the second follows from
Proposition 5.1 of Kl\"{u}ppelberg et al. (2004).

From Kl\"{u}ppelberg et al.  (2004, Theorem 4.2), we can deduce:

\begin{theorem}\label{ceth1}
Assume $\alpha> 0$ and \eqref{c1} holds.
Then for $x\ge 0$
\begin{equation}\label{cet2a}
\lim_{u \to \infty} P^{(u)}\left(X_{\tau_u}-u\le x\right)=
1-\beta_2 e^{-\alpha x}
-\frac{1}{q}\int_{(0,\infty)}\left(e^{\alpha y}-1\right)
\Pi_{\cal H}(x+\rmd y).
\end{equation}
\end{theorem}
\medskip

Equation \eqref{cet2a} compares closely with the result in \eqref{t2a}, and reduces to it if we formally set $\alpha=\nu_0$
and  $Ee^{{\alpha}{X_1}}=1$, since in that case  $\beta_2=0$ by \eqref{did}.
For $x=0$, using \eqref{kapexp} and \eqref{betadef}, \eqref{cet2a} gives
\begin{equation}\label{acr}
\lim_{u \to \infty} P^{(u)}\left(X_{\tau_u}=u\right)
=\frac{\alpha\rmd_H}{q}
\end{equation}
as the asymptotic (conditional) probability of creeping in the
convolution equivalent case.
Compare this with the analogous result in the Cram\'er case given by \eqref{t3a}.
If $\rmd_H>0$, then  \eqref{kes} and \eqref{acr} imply
\begin{equation}\label{advlim}
\lim_{u \to \infty}
\frac{V'_H(u)}{P(\tau_u<\infty)}
=\frac{\alpha}{q},
\end{equation}
with the analogous result, with $\alpha$ replaced by $\nu_0$, following from \eqref{t3a} in the Cram\'er case.

The limiting distributions of the undershoots in the convolution equivalent case are given in Table 1 below.
Their derivations using the quintuple law are explained in Doney and Kyprianou (2006) and carried out
in Park and Maller (2008).
An alternative approach, yielding more general results, is given in Griffin and Maller (2010).

\setcounter{equation}{0}\section{Summary of General Results}\label{s6}
The general results are summarised in Table 1. In it, ``Cram\'er Case"
refers to results proved under Assumptions \eqref{cramer}
and \eqref{mu*},
while the ``Convolution Equivalent Case"
refers to results proved under Assumption \eqref{c1}, with
$\alpha>0$. All the asymptotic distributions are valid for $x\ge 0$. Recall \eqref{betadef},
and note that $\beta_2=0$ in the
 Cram\'er case.
To the entries in the table we can also add:

\bigskip \noindent {\bf Cram\'er Case:}
\begin{equation}\label{rCr}
P(\tau_u<\infty)
\sim  \frac{qe^{-\nu_0 u}}{\nu_0 m^*}, \ {\rm as}\ u\to\infty.
 \end{equation}

\medskip\noindent  {\bf Convolution Equivalent Case:}
\begin{equation}\label{rCE}
P(\tau_u<\infty)
\sim
\frac{1}{\beta_1\beta_2}
\overline{\Pi}_{X}^+(u)
\sim
\frac{1}{q\beta_2^2}
\overline{\Pi}_{\cal H}(u), \ {\rm as}\ u\to\infty.
\end{equation}
The result in the convolution equivalent
 case follows from Theorem 4.1 and Proposition 5.3 in  Kl\"{u}ppelberg et al. (2004).
Note that in the Cram\'er case, $\overline{\Pi}_{X}^+(u)$ and  $\overline{\Pi}_{\cal H}(u)$ are both of
smaller order than $e^{-\nu_0 u}$ as $u\to\infty$,
by \eqref{mu*} and \eqref{picon}.

\begin{table}[h!]
\begin{narrow}{-2cm}{-2cm}
\begin{center}
    \footnotesize{
    \begin{tabular}{| c | c | c | c |}
    \hline
    & Limiting distribution              & Case I:                     & Case II:                  \\
    & Valid for all $x \geq 0$           & Cram{\'e}r case             & Convolution Equivalent Case  \\
    \hline
  I \KeesT \KeesB & $\displaystyle \lim_{u \to \infty} P^{(u)}(X_{\tau_u}-u \leq x)$   & $\displaystyle 1-\dfrac{1}{q} \int_{(0,\infty)} (e^{\nu_0 y}-1) \Pi_{\mathcal{H}}(x+\mathrm{d}y)$ & $\displaystyle 1-\beta_2 e^{-\alpha x}-\dfrac{1}{q} \int_{(0,\infty)} (e^{\alpha y}-1) \Pi_{\mathcal{H}}(x+\mathrm{d}y)$ \\ 
    \hline
  II \KeesT \KeesB & $\displaystyle \lim_{u \to \infty} P^{(u)}(u-X_{\tau_u-} \leq x)$   & $\displaystyle \frac{\nu_0 \mathrm{d}_{H}}{q} + \frac{1}{q} \int_0^x g_x(y) \, \mathrm{d}(e^{\nu_0 y})$ & $\displaystyle \frac{\alpha \mathrm{d}_{H}}{q} + \frac{1}{q} \int_0^x g_x(y) \, \mathrm{d}(e^{\alpha y})$ \\ 
    \hline
  III \KeesT \KeesB & $\displaystyle \lim_{u \to \infty} P^{(u)}(u-\overline{X}_{\tau_u-} \leq x)$   & $\displaystyle \frac{\nu_0 \mathrm{d}_{H}}{q} + \frac{1}{q} \int_0^x \overline{\Pi}_{\mathcal{H}}(y) \, \mathrm{d}(e^{\nu_0 y})$ & $\displaystyle \frac{\alpha \mathrm{d}_{H}}{q} + \frac{1}{q} \int_0^x \overline{\Pi}_{\mathcal{H}}(y) \, \mathrm{d}(e^{\alpha y})$ \\ 
    \hline
  IV \KeesT & \KeesB $\displaystyle \lim_{u \to \infty} P^{(u)}(X_{\tau_u} =u)$   & $\displaystyle \frac{\nu_0 \mathrm{d}_H}{q}$ & $\displaystyle \frac{\alpha \mathrm{d}_H}{q}$ \\ 
    \hline
    \end{tabular}
    }
\end{center}
\end{narrow}
  \caption{Limiting Distributions of Overshoot and Undershoots in the Cram\'er and Convolution Equivalent Cases}
\end{table}

\medskip\noindent
{\bf Remark:}\
Table 1 corrects a couple of expressions in Park and Maller (2008).
The proof of their Theorem 3.2 is correct but there are errors in their Eq. (3.5);
the simpler expression in Row II Column 3 of Table 1 should be used instead.
Also, Eq. (3.10) of Park and Maller (2008) should be corrected to that shown in
Row III Column 3 of Table 1 by adding
in the term from creeping.
\medskip

It is important to note that in the  convolution equivalent case,
the limiting distributions of the undershoots are improper.
Letting $x\to\infty$ and using \eqref{ginf}, we see that they both have total mass given by
\ben\ba
\frac{\alpha\rmd_H}{q}
+ \frac{1}{q}\int_0^\infty
\pibar_{\cal H}(y)\rmd(e^{\alpha y})
&=\frac{\alpha\rmd_H}{q}
+ \frac{1}{q}\int_0^\infty (e^{\alpha y}-1)
\pibar_{\cal H}(\rmd y)\\
&=\frac {q-\kappa(0,-\alpha)}{q}=1-\beta_2<1.
\ea\een
The remaining mass escapes to $\infty$ because with positive probability, as $u\to\infty$,
$X$ can pass over level $u$ as the result of an exceedingly  large jump
the first time it leaves a neighbourhood of the origin.
See Doney and Kyprianou (2006) and Griffin and Maller (2010) for a more detailed discussion.

\section{Example: the GTSC-Class}\label{s7}
\setcounter{equation}{0}
In this section we use the formulae listed in Table 1 to examine the asymptotic distributions of the
 overshoot and undershoots,
for a particular class, the GTSC (``Gaussian-Tempered-Stable-Convolution") class, of L\'evy processes,
introduced by Hubalek and Kyprianou (2010).
For these, the marginal distributions of $X_1$ are of generalised inverse Gaussian type.
They generate a class of L\'evy processes
$(X_t)_{t\ge 0}$ which are spectrally positive with corresponding upward ladder height processes,
$({{H}}_t)_{t\ge 0}$, being killed tempered stable subordinators.
The downward ladder height processes are trivial subordinators normalised to unit drift.
These processes admit some quite explicit expressions for the
L\'evy measures
of $X$ and of ${{H}}$, and so are very suitable for calculation purposes.
Our main aim is to illustrate similarities and differences between asymptotic distributions
of the overshoot and undershoots for the Cram\'er and convolution equivalent cases.

The GTSC class depends on five parameters
$q>0$, $\rmd_{{H}}\ge 0$, $c>0$, $\alpha \ge 0$ and $\rho\in [-1,1)$.
When $\rho \leq 0$  we  require
$\alpha>0$.
The Laplace exponent of $X$ takes the following forms:
\begin{equation}\label{def_phi_X}
\psi_X(\theta) = E[e^{\theta X_1}]=
-q \theta+\rmd_{{H}}\theta^2
-c\theta\Gamma(-\rho)\left(\alpha^\rho-(\alpha-\theta)^\rho \right),\
\mbox{for $\theta<\alpha$},\ \rho\ne 0
\end{equation}
and
\begin{equation}\label{1.21}
\psi_X(\theta) = -q\theta + \rmd_{{H}}\theta^2
+ c\theta\log \frac{\alpha}{\alpha-\theta},\  \mbox{for $\theta<\alpha$},\ \rho= 0.
\end{equation}
It can be checked that the L\'evy measure of $X$ is given by
\begin{equation}\label{def_Pi_X}
\Pi_X(\rmd y) = c\left( \frac{\alpha}{y^{\rho+1}}+\frac{\rho+1}{y^{\rho+2}}\right)e^{-\alpha y}\rmd y, \quad y>0,
\end{equation}
hence $\overline{\Pi}_{X}^+ \in {\cal S}^{(\alpha)}$ if
$\alpha>0$ and $\rho>0$, by the result of Chover et al. (1973) referred to earlier.
We use $\alpha$ for the fourth parameter of the GTSC class for this reason,
so as to connect with the results in the previous sections.
In terms of modelling, the parameters $c$, $\alpha$ and $\rho$ determine the claim size distribution, $q$ the premium rate and $\rmd_H$ the noise or Brownian perturbation.

The ascending ladder height process ${H}$ is a subordinator
killed at rate $q$, with drift $\rmd_{{H}}$ and L\'evy measure
\begin{equation}\label{1.22}
\Pi_{{\cal H}}(\rmd y) = c y^{-\rho-1} e^{-\alpha y} \, \rmd y, \quad y>0.
\end{equation}
We have $\overline{\Pi}_{\cal H} \in {\cal S}^{(\alpha)}$ if, again,  $\alpha>0$ and $\rho>0$.
The Laplace exponent of $H$ is given by
\begin{equation}\ba\label{def_phi_H}
\psi_{{H}}(\theta)
&=-q+\rmd_{{H}}\theta-c\Gamma(-\rho)\left(\alpha^{\rho}-(\alpha-\theta)^{\rho}\right),\  \mbox{for $\theta<\alpha$},\ \rho\ne 0
\ea\end{equation}
and
\begin{equation}\label{def_phi_H0}
\psi_{{H}}(\alpha)=
-q + \rmd_{{H}}\theta
+ c\log \frac{\alpha}{\alpha-\theta},\  \mbox{for $\theta<\alpha$},\ \rho= 0.
\end{equation}
For $\rho \in [0,1)$, $H$ has infinite activity ($\Pi_{\cal H}(\R)=\infty$).
For $\rho=0$ it is simply a (killed) Gamma subordinator (with drift).
If $-1\le \rho<0$ it is a (killed) compound Poisson process (with drift),
with intensity parameter $c\alpha \rho\Gamma(-\rho)$ and
Gamma distributed jumps.
For details of all these properties see Hubalek and Kyprianou (2010).

From \eqref{def_phi_X} and \eqref{1.21} we see that
\[
E[X_1]=\psi_X'(0)=-q<0,
\]
so $X_t \to -\infty$ a.s. as $t \to \infty$.
Recall that the Cram\'er case is characterised by \eqref{cramer} and \eqref{mu*},
while the convolution equivalent case is characterised by \eqref{c1}.
For particular combinations of the GTSC parameters, the Cram\'er case is also
included in the GTSC class.  We can classify $X$ as follows.

\begin{theorem}\label{Keesth1}
For $\alpha>0$, and $X$ from the GTSC-class parameterised as in \eqref{def_phi_X} or \eqref{1.21},
we have the following:
\begin{itemize}
\item[(i)] if either $\rho \in [-1,0]$, or $\rho \in (0,1)$ and
$\rmd_{{H}} \alpha -q -c\alpha^{\rho} \Gamma(-\rho)>0$,
$X$ is in  the Cram\'er case  with parameter $\nu_0$,
where $\nu_0$ is the unique root of $\psi_X=0$ on the interval $(0,\alpha)$;
thus, $\nu_0$ satisfies
\begin{equation}\label{root}
0=\psi_X(\nu_0)/\nu_0=
\rmd_{{H}}\nu_0-q-c\Gamma(-\rho)\left(\alpha^{\rho}-(\alpha-\nu_0)^{\rho}\right),\
0<\nu_0<\alpha;
\end{equation}
\item[(ii)] if $\rho \in (0,1)$ and $\rmd_{{H}} \alpha -q-c\alpha^{\rho} \Gamma(-\rho)<0$,
$X$ is in  the convolution equivalent  case  with parameter $\alpha$;
\item[(iii)] if $\rho \in (0,1)$ and $\rmd_{{H}} \alpha -q-c\alpha^{\rho} \Gamma(-\rho)=0$,
$X$ is in neither case.
\end{itemize}
\end{theorem}

\bigskip \noindent {\bf Proof of Theorem \ref{Keesth1}}.\
(i)\
Recall that
\be\label{psiprop}
\text{$\psi_X$ is convex, $\psi_X(0)=0$ and $\psi_X'(0)<0$.}
\ee
Hence, if $\psi_X(\alpha-)>0$ it follows that
$\psi_X$ has a unique zero $\nu_0$ on $(0,\alpha)$ and for $\veps>0$ small enough also
$\psi_X(\nu_0+\veps)<\infty$, giving \eqref{cramer} and \eqref{mu*}.
It is straightforward to check that the conditions in (i) indeed yield $\psi_X(\alpha-)>0$.

(ii)\
Next, assume $\rho \in (0,1)$ and $\rmd_{{H}} \alpha -  q-c\alpha^{\rho} \Gamma(-\rho)<0$.
Then
$\overline{\Pi}_{X}^+ \in {\cal S}^{(\alpha)}$  since $\rho>0$, and further
\[
\log E[e^{\alpha {X}_1}]=
\psi_{{X}}(\alpha-)
=\alpha\left(\rmd_{{H}} \alpha-q-c\alpha^{\rho} \Gamma(-\rho)\right)<0.
\]
Thus \eqref{c1} holds.

(iii)\
Finally, suppose $\rho \in (0,1)$ and $\rmd_{{H}} \alpha -  q -c\alpha^{\rho} \Gamma(-\rho)=0$.
We are not in  the convolution equivalent  case  since by the same argument as in (ii) we now have
$\log E[e^{\alpha {X_1}}]= 0$. Nor are we in  the Cram\'er case,
since the only candidate for $\nu_0$ is $\alpha$, by \eqref{psiprop}. But then it follows immediately from \eqref{rid} that $m^*=\infty$ and hence \eqref{mu*} fails by \eqref{mu*1}.
\hfill\halmos\bigskip

It will suffice for our demonstrations {\em to keep $\alpha >0$ and $\rho\in(0,1)$ in what follows.}
Thus we restrict to the parameter set
$q>0$, $\rmd_{{H}}\ge 0$, $c>0$, $\alpha> 0$ and $\rho\in (0,1)$.
The parameter $\beta_2$ defined by \eqref{betadef}  here takes the form:
\begin{equation}\label{10nov5}
\beta_2= \frac{ \kappa(0,-\alpha)}{ q }
= \frac{ -\psi_{{H}}(\alpha)}{ q }
= \frac{ q -\alpha \rmd_{{H}} +c \alpha^\rho\Gamma(-\rho)}{ q }.
\end{equation}
Thus by Theorem \ref{Keesth1},  the Cram\'er case arises when $\beta_2<0$ and the convolution equivalent case
when $\beta_2>0$.  The boundary between these regions is the surface defined by
\be\ba
B_0:=\{p=(q, \rmd_{{H}}, c, \alpha, \rho): \beta_2=0\}.
\ea\ee
We will denote the regions in parameter space where $\beta_2<0$
and $\beta_2>0$ by $B_-$ and $B_+$ respectively.

We will consider the asymptotics relevant to the GTSC class,
using phrases like ``approach the boundary through the
Cram\'er class (convolution equivalent class)"
and variations thereof to mean that we
continuously move through points $p\in B_-$ ($p\in B_+$) to a point $p_0\in B_0$.
The following observations will be used in the discussion below without further mention.
If $p^0=(q^0, \rmd_{{H}}^0, c^0, \alpha^0, \rho^0)\in B_0$ and we approach $p^0$ through points
$p=(q, \rmd_{{H}}, c, \alpha, \rho)\in B_-$, then one easily checks, using \eqref{root},
that $\nu_0\to \alpha^0$. Hence $m^* \to \infty$ by \eqref{rid}.
Of course, trivially, $\alpha\to\alpha^0$ and  $\beta_2\to 0$ as $p\to p^0$
irrespective of whether $p\in B_-$ or $p\in B_+$.

The following identity, which follows easily from the definition of the gamma function, is useful in deriving some of the formulae given in the remainder of this section;
for $0<\theta\le \alpha$ and $\rho\in (0,1)$
\begin{equation}\label{intid}
\int_0^\infty (e^{\theta y}-1) y^{-\rho-1} e^{-\alpha y}\rmd y
=-\Gamma(-\rho)\left(\alpha^{\rho}-(\alpha-\theta)^{\rho}\right).
\end{equation}

For the probability of ruin
we use the formulae in \eqref{rCr} and \eqref{rCE}.
In  the Cram\'er case, by \eqref{rid} and  \eqref{1.22}, we have
\begin{equation}\label{10nov1}
P\left(\tau_u<\infty\right)
\sim \frac{q(\alpha-\nu_0)^{1-\rho}}{\nu_0\left(\rmd_H(\alpha-\nu_0)^{1-\rho}-c\rho\Gamma(-\rho)\right)}e^{-\nu_0 u}, \quad \mbox{as $u \to \infty$}
\end{equation}
while in  the convolution equivalent case, by \eqref{10nov5}, we have
\begin{equation}\label{10nov2}
P\left(\tau_u<\infty\right)
\sim \frac{cq}{\alpha\left( q-\alpha\rmd_H+c\alpha^\rho\Gamma(-\rho)\right)^2}u^{-\rho-1}e^{-\alpha u},  \quad
\mbox{as $u \to \infty$.}
\end{equation}
There are clear structural differences between these estimates.
For example, in  the Cram\'er case,  the rate of decay of the ruin probability depends,
through $\nu_0 \in (0,\alpha)$,  on all of the parameters $q$, $\rmd_H$, $c$, $\alpha$ and $\rho$,
while in  the convolution equivalent  case  it depends only on $\alpha$ and $\rho$.

Now fix $p^0=(q^0, \rmd_{{H}}^0, c^0, \alpha^0, \rho^0)\in B_0$ and  approach $p^0$ through points
 $p\in B_-$ or points $p\in B_+$. If $p\in B_-$, then $\nu_0\to \alpha^0$,
while if $p\in B_+$, then $\alpha\to \alpha^0$.
Thus the rate of  exponential decay
in \eqref{10nov1} and \eqref{10nov2}
changes continuously across the boundary.
The constant in (\ref{10nov1}) vanishes as we approach the boundary since $m^* \to \infty$,
while the constant in  \eqref{10nov2} explodes  since $\beta_2\to 0$.
This, of course,  is a consequence of the structurally different types of decay in the two cases.

The formulae for the asymptotic distributions of the overshoot and undershoots under the GTSC assumption are listed in Table 2.  These are
calculated using the corresponding quantities in Table 1.

\begin{table}[h!]
\begin{narrow}{-2cm}{-2cm}
\begin{center}
    \footnotesize{
    \begin{tabular}{| c | c | c | c |}
    \hline
    & Limiting distribution              & Case I:                     & Case II:                  \\
    & Valid for all $x \geq 0$           & Cram{\'e}r case             & Convolution Equivalent Case  \\
    \hline
  I \KeesT \KeesB & $\displaystyle \lim_{u \to \infty} P^{(u)}(X_{\tau_u}-u \leq x)$   & $\displaystyle 1-\frac{c}{q} \int_x^{\infty} y^{-\rho-1} e^{-\alpha y} \left( e^{\nu_0 (y-x)}-1 \right) \, \mathrm{d}y$ & $\displaystyle 1-\beta_2 e^{-\alpha x}-\frac{c}{q} \int_x^{\infty} y^{-\rho-1} \left( e^{-\alpha x}-e^{-\alpha y} \right) \, \mathrm{d}y$ \\ 
    \hline
  II \KeesT \KeesB & $\displaystyle \lim_{u \to \infty} P^{(u)}(u-X_{\tau_u-} \leq x)$   & $\displaystyle \frac{\nu_0 \mathrm{d}_{H}}{q} + \frac{c}{q} \int_0^x 
y^{-\rho-1} e^{-\alpha y} \left( e^{\nu_0 y}-1 \right) \, \mathrm{d}y$ & $\displaystyle \frac{\alpha \mathrm{d}_{H}}{q} + \frac{c}{q} \int_0^x y^{-\rho-1} \left( 1-e^{-\alpha y} \right) \, \mathrm{d}y$ \\ 
    \hline
  III \KeesT \KeesB & $\displaystyle \lim_{u \to \infty} P^{(u)}(u-\overline{X}_{\tau_u-} \leq x)$   & $\displaystyle 1-\frac{c}{q} \int_x^{\infty} y^{-\rho-1} e^{-\alpha y} \left( e^{\nu_0 y}-e^{\nu_0 x} \right) \, \mathrm{d}y$ & $\displaystyle 1-\beta_2-\frac{c}{q} \int_x^{\infty} y^{-\rho-1} e^{-\alpha y} \left( e^{\alpha y}-e^{\alpha x} \right) \, \mathrm{d}y$ \\ 
    \hline
  IV \KeesT & \KeesB $\displaystyle \lim_{u \to \infty} P^{(u)}(X_{\tau_u} =u)$   & $\displaystyle 1+\frac{c}{q} \Gamma(-\rho) \left( \alpha^{\rho} -(\alpha-\nu_0)^{\rho} \right)$ & $\displaystyle 1-\beta_2+\frac{c \alpha^{\rho}}{q} \Gamma(-\rho)$ \\ 
    \hline
    \end{tabular}
    }
\end{center}
\end{narrow}
  \caption{Limiting Distributions of Overshoot and Undershoots for the GTSC Class}
\end{table}

Theorem \ref{Keesth1} shows that the sign of the quantity $\rmd_{{H}} \alpha -  q -c\alpha^{\rho} \Gamma(-\rho)$
distinguishes between the Cram\'er and convolution equivalent cases.
In Fig.1 we plot it as a function of $\alpha$, for a particular choice of parameters
(which are used for all 4 figures); namely, fix
\begin{equation}\label{1}
q=1,   {\rmd}_{{H}}=1/2, c=1, \rho=1/2,
\end{equation}
and consider the function
\begin{equation}\label{fdef}
f(\alpha):=
\rmd_{{H}} \alpha -  q -c\alpha^{\rho} \Gamma(-\rho)
\end{equation}
for $\alpha \in (0,1)$.
{}From Fig. 1 we see that $f$ is negative on an interval $(0,\alpha_0)$,
the Cram\'er case,  and positive on $(\alpha_0,1)$, the convolution equivalent case,
where, for the parameter values in (\ref{1}), $\alpha_0=0.069$.

\begin{figure}[!ht]
  \centering
    \includegraphics[width=0.7\textwidth]{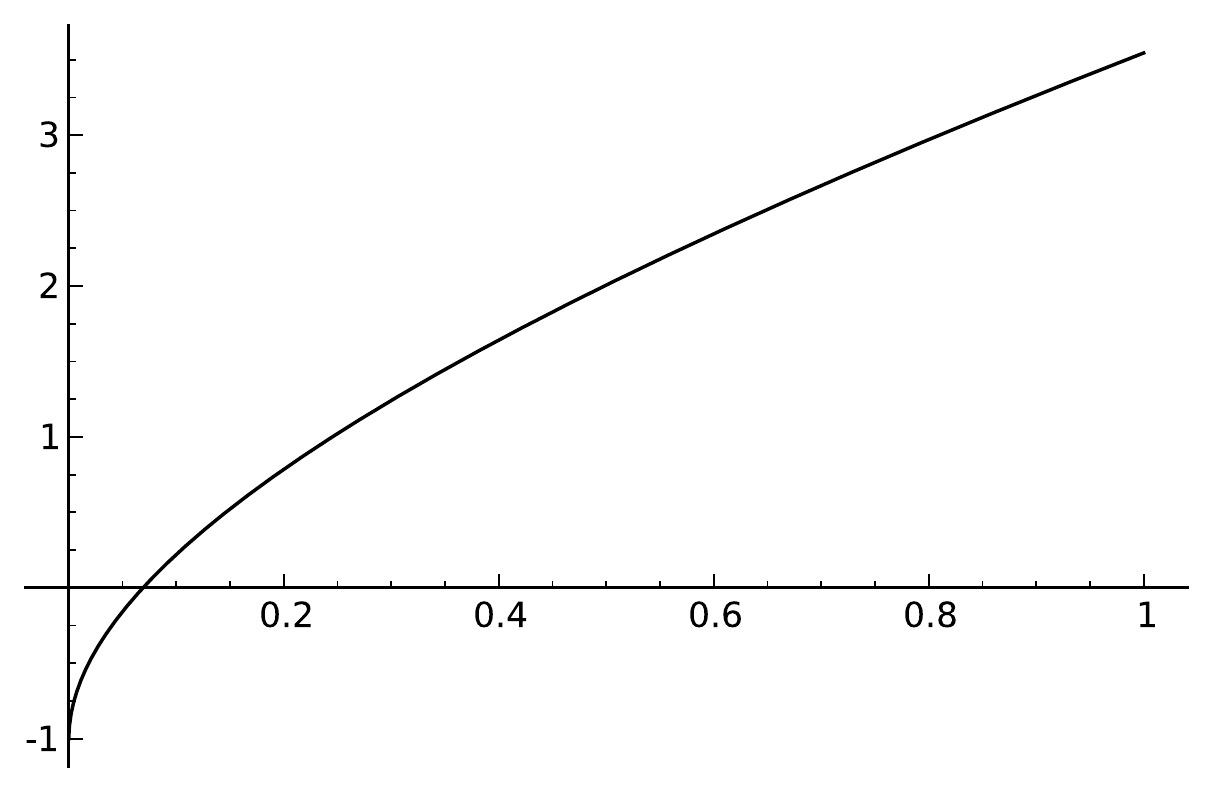}
\caption{The function $f$ defined in (\ref{fdef}) for values $0<\alpha<1$.
Parameter values are $q=1$, $\rmd_H=1/2$, $c=1$, $\rho=1/2$.}
\end{figure}

The cumulative distribution functions of the overshoot,  for both cases,
Cram\'er and convolution equivalent, are represented in the next diagram as $\alpha$ varies,
for the same choice of parameters as in Fig 1.
Recall that varying $\alpha$ is a means of varying the claim size distribution.  Even though increasing $\alpha$ beyond $\alpha_0$ results in changing the model from one in which the Cram\'er condition holds to one in which it doesn't, the overshoot distribution transitions smoothly throughout the entire range of $\alpha$.  A similar phenomenon can be observed by varying the other parameters.

\begin{figure}[!ht]
  \centering
    \includegraphics[width=0.7\textwidth]{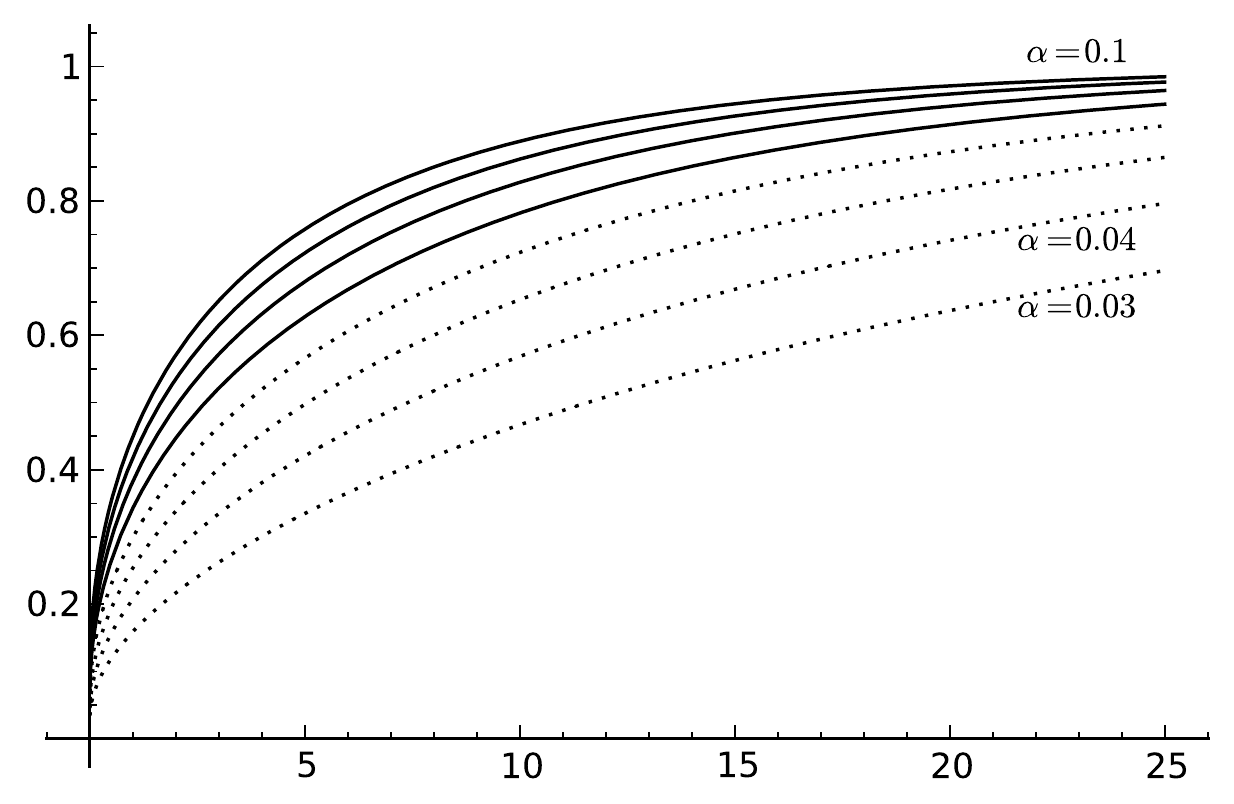}
\caption{Cumulative distribution of the overshoot as given by Row I of Table 2,
 in the Cram\'er case for values of $\alpha=0.07, 0.08, 0.09, 0.10  $ (full lines),
and in the convolution equivalent case for values of $\alpha=0.03, 0.04, 0.05, 0.06$ (dotted lines).
Parameter values as for Fig. 1.}
\end{figure}

The formulae for the overshoot in Row I of Table 2 are obtained by direct substitution of
\eqref{1.22} into Table 1.
For computational purposes these formulae can be expressed in terms of the
incomplete Gamma function, and easily calculated
for specific values of the parameters.

As the boundary is approached, either through $B_-$ or $B_+$,
the formulae agree in the limit, for every $x\ge 0$.  Figure 2 illustrates this in the case that the boundary is approached by varying $\alpha$ only.
Setting $x =0$ gives the asymptotic (conditional) probability of creeping over the level $u$.
The resulting formulae agree with those in Row IV of Table 1:
In the Cram\'er case
\begin{eqnarray}\label{10nov6}
\lim_{u \to \infty} P^{(u)}\left(X_{\tau_u}=u\right)
&=&
 1 - \frac{c}{ q } \int_0^\infty (e^{\nu_0 y}-1) y^{-\rho-1} e^{-\alpha y}\rmd y
 \nonumber\\
&=&
1+\frac{c}{ q } \Gamma(-\rho)\left(\alpha^{\rho}-(\alpha-\nu_0)^{\rho}\right)=\frac{\rmd_{{H}} \nu_0}{q},
\end{eqnarray}
where the second equality holds by \eqref{intid} and the third by \eqref{root}, while
for  the convolution equivalent  case, by the same means but using (\ref{10nov5}), we find that
\begin{equation}\label{10nov7}
\lim_{u \to \infty} P^{(u)}\left(X_{\tau_u}=u\right) =
1 - \beta_2 +\frac{c\alpha^{\rho}}{ q } \Gamma(-\rho)
= \frac{\rmd_{{H}} \alpha}{ q }.
\end{equation}

Comparing (\ref{10nov6}) and (\ref{10nov7}), note that in both cases there is creeping in the limit if and
only if $\rmd_{{H}}>0$, which is in turn equivalent to $X$ having a Gaussian component
(namely, from (\ref{def_phi_X}) we see that $\sigma_X=\sqrt{2\rmd_{{H}}}$).
One obvious  difference between the two cases  is that in
the Cram\'er case  changing any parameter changes $\nu_0$ and hence (\ref{10nov6}),
while (\ref{10nov7}) does not depend on $\rho$ and $c$ (we saw a similar effect for the ruin probability).

Next, for the behaviour when $x \to \infty$,
we can calculate, in  the Cram\'er case
\[
1-\lim_{u \to \infty} P^{(u)}\left(X_{\tau_u}-u \leq x\right)
\sim
\frac{c}{ q }\left((\alpha-\nu_0)^{-1}-\alpha^{-1}\right)
x^{-\rho-1} e^{-\alpha x},\ x \to \infty,
 \]
while in  the convolution equivalent  case
\[
1-\lim_{u \to \infty} P^{(u)}\left(X_{\tau_u}-u \leq x\right)
\sim \beta_2 e^{-\alpha x} \quad \mbox{as $x \to \infty$}.
\]
The structural difference between the two cases is largely analogous to that seen for the ruin probabilities,
but curiously with the roles reversed.
In this case, approaching the boundary through the Cram\'er class causes the coefficient to explode,
whereas the coefficient vanishes if the approach is through the convolution equivalent class.

This analysis shows that,
while the formulae in the two cases are structurally different,
especially in their asymptotic behaviour, nevertheless they segue continuously into each other across the
boundary $B_0$. This effect is apparent in Fig 2.

Similar analyses can be given for the undershoots.

\begin{figure}[!ht]
  \centering
    \includegraphics[width=0.7\textwidth]{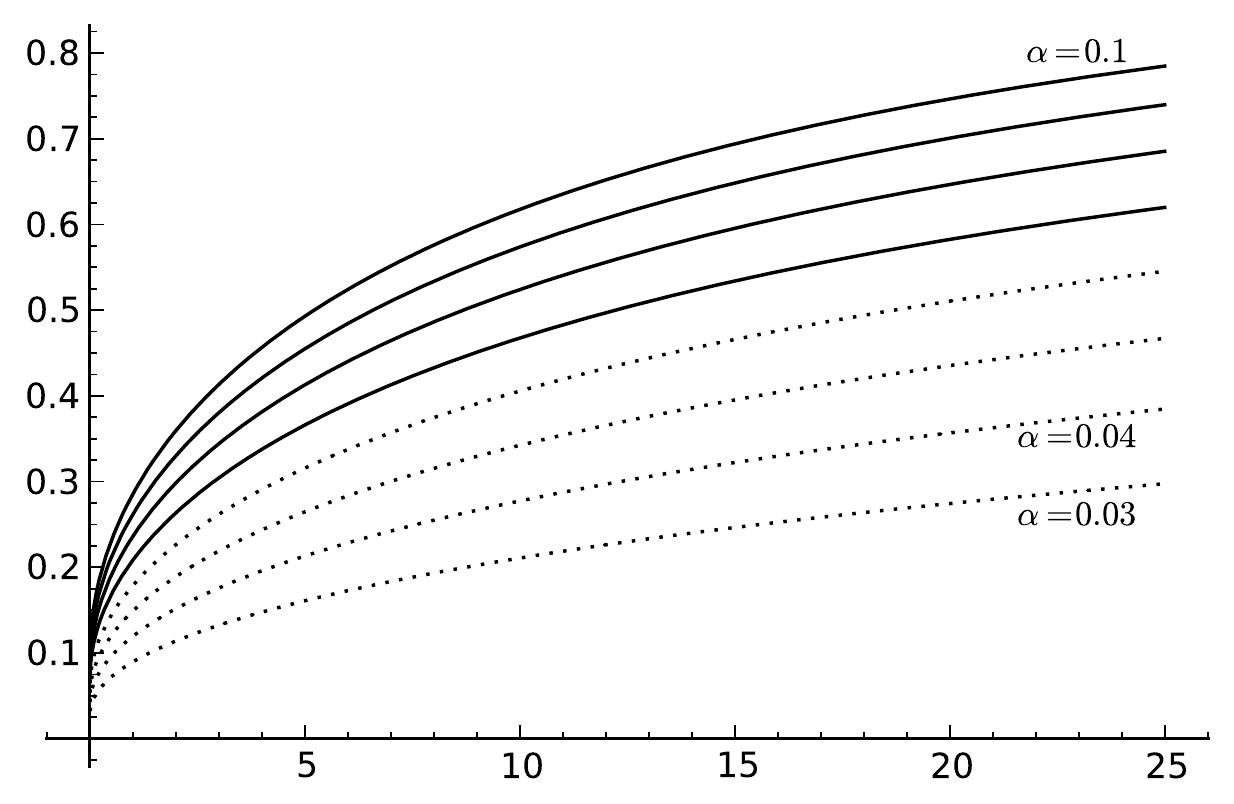}
\caption{Cumulative distribution of the undershoot as given by Row II of Table 2.
Values of $\alpha$ and parameter values as for Fig. 1.}
\end{figure}

For the undershoot itself, from \eqref{DKcor6}, \eqref{gdef} and the fact that $\wh V_{\wh H}(\rmd y)=\rmd y$ since $\wh H$ is unit drift, we obtain
\ben
g_x(\rmd y)=-\pibar_X^+(y)\mbox{d}y=-\Pi_{\cal H}(\mbox{d}y).
\een
Integrating by parts in Row II of Table 1, then yields
the results in Row II of Table 2:
\begin{equation}\label{17nov1}
\lim_{u \to \infty} P^{(u)}(u-X_{\tau_u-} \leq x) = \frac{\nu_0 \rmd_{\mathcal{H}}}{q}
+\frac{c}{ q }\int_0^x y^{-\rho-1} e^{-\alpha y} \left( e^{\nu_0 y}-1 \right)\mbox{d}y
\end{equation}
in the Cram\'er case, and
\begin{equation}\label{17nov2}
\lim_{u \to \infty} P^{(u)}\left(u-X_{\tau_u-} \leq x\right)
= \frac{\alpha \rmd_{\mathcal{H}}}{q}
+\frac{c}{ q }\int_0^x y^{-\rho-1} \left(1- e^{-\alpha y}\right)\mbox{d}y
\end{equation}
in the convolution equivalent case.
Just as for the overshoot, there is a continuous transition across the boundary in the formulae
for every $x\ge 0$.  When $x=0$ they reduce to
\begin{eqnarray}\label{17nov5c}
\lim_{u \to \infty} P^{(u)}\left({X}_{\tau_u-}=u\right)
=\frac{\nu_0 \rmd_{{H}}}{q}
\ \ \text{ and }\ \
\lim_{u \to \infty}P^{(u)}\left({X}_{\tau_u-}=u\right)
=\frac{\alpha \rmd_{{H}}}{q}
\end{eqnarray}
respectively.
When  $x \to \infty$,  the right hand side of (\ref{17nov1}) converges to
\begin{equation}\label{721}
\frac{\nu_0 \rmd_{{H}}}{q}-\frac{c}{q}\Gamma(-\rho)\left(\alpha^{\rho}-(\alpha-\nu_0)^{\rho}\right) =1,
\end{equation}
by  \eqref{root} and \eqref{intid},
while the right hand side of (\ref{17nov2}) converges to
\begin{equation}\label{722}
\frac{\alpha \rmd_{{H}}}{ q }-\frac{c\alpha^{\rho}\Gamma(-\rho)}{ q }
=1-\beta_2
\end{equation}
by \eqref{10nov5} and \eqref{intid}.  Since $\beta_2>0$ in the convolution equivalent case,
this last expression is in $(0,1)$.
Hence in the convolution equivalent case, the undershoot converges under
$P^{(u)}$ to an improper distribution, which has an atom at $\infty$ of mass $\beta_2$.
This atom vanishes as we approach the boundary through $B_+$.

Just as for the overshoot, another difference between the cases is in their asymptotic behaviour, namely, in the Cram\'er case we see
\be\label{4J1}
1-\lim_{u \to \infty} P^{(u)}\left(u-X_{\tau_u-}\leq x\right)
\sim \frac{c }{ q (\alpha-\nu_0)} x^{-\rho-1} e^{-(\alpha-\nu_0)x},
\quad \mbox{as $x \to \infty$},
 \ee
while in the convolution equivalent case
\be\label{4J2}
1-\beta_2
-\lim_{u \to \infty} P^{(u)}\left(u-X_{\tau_u-}\leq x\right)
\sim \frac{c}{\rho  q } x^{-\rho}, \quad \mbox{as $x \to \infty$.}
\ee
As the boundary is approached through $B_-$, the exponential factor in \eqref{4J1} vanishes, but the coefficient explodes.  In contrast, the estimate in \eqref{4J2} behaves smoothly at the boundary.

For the undershoot from the previous maximum,
the formulae in Row III of Table 2  follow from Row III of Table 1
after an integration by parts and substitution of \eqref{1.22}.
Just as for the overshoot and the undershoot, there is a continuous transition across the boundary in the formulae  for every $x\ge 0$.  When $x=0$, upon using \eqref{intid}, they reduce to
\begin{eqnarray}
\lim_{u \to \infty} P^{(u)}\left(\overline{X}_{\tau_u-}=u\right)
=\frac{\nu_0 \rmd_{{H}}}{q}
\ \ \text{and}\ \
\lim_{u \to \infty}P^{(u)}\left(\overline{X}_{\tau_u-}=u\right)
=\frac{\alpha \rmd_{{H}}}{q}
\end{eqnarray}
respectively.
Precisely the same asymptotics as $x \to \infty$ apply as in \eqref{721} and \eqref{722}.
  As for the undershoot,
in the convolution equivalent case, the undershoot of the last maximum before $\tau_u$ converges under
$P^{(u)}$ to an improper distribution, with an atom at $\infty$ of mass $\beta_2$.
This atom vanishes as we approach the boundary through $B_+$.  Finally \eqref{4J2} continues to hold in the convolution equivalent case if ${X}_{\tau_u-}$ is replaced by $\overline{X}_{\tau_u-}$, while in the Cram\'er case
\ben
1-\lim_{u \to \infty} P^{(u)}\left(u-\overline{X}_{\tau_u-}\leq x\right)
\sim \frac{c }{ q} ((\alpha-\nu_0)^{-1}-\alpha^{-1}) x^{-\rho-1} e^{-(\alpha-\nu_0)x},
\quad \mbox{as $x \to \infty$}.
 \een

\begin{figure}[!ht]
  \centering
    \includegraphics[width=0.7\textwidth]{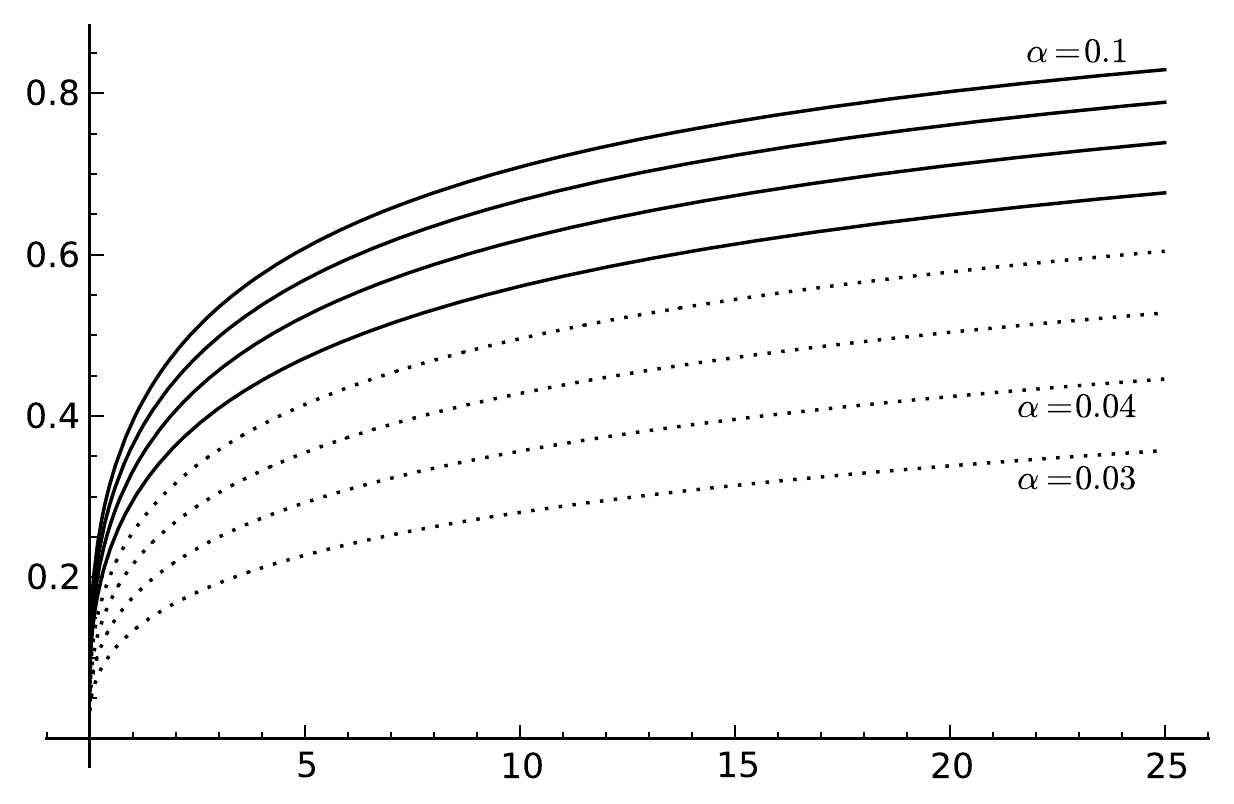}
\caption{Cumulative distribution of the undershoot from the previous maximum as given by Row III of Table 2.
Values of $\alpha$ and parameter values as for Fig. 1.}
\end{figure}

\setcounter{equation}{0}\section{Discussion}\label{s8}
The quintuple law, referenced in Section \ref{s3}, is the result of  a deep analysis into L\'evy process theory,
but the derivations from it of distributional identities such as those in Section \ref{s3}, and, from them, the applications to limit laws such as those  in
Sections \ref{s4} and  \ref{s5},  are straightforward.  In fact, once given the results in Section \ref{s3},
our paper is virtually self-contained, requiring in addition nothing more than some standard renewal theory and
the dominated convergence theorem for the proof of Theorem \ref{th1}.
The proof of Theorem \ref{ceth1}, also based on the quintuple law (see Kl\"uppelberg et al. (2004)),
requires some extra techniques from the theory of convolution equivalent distributions,
but these are standard and easily applied.

Consequently, we can put forward the methods exemplified in this paper as a simple and clear way
to set out and derive limit laws for overshoots and undershoots in both the convolution equivalent and
Cram\'er formulations. On the other hand, the method is not a panacea for all problems.
Inspection of Eq. \eqref{cdf4} for the joint distribution of the 3 positions variables, and of
Eq. \eqref{cdf6a} for the joint distribution of the 2 time variables, suggests that our methods will
not carry over easily to these situations. We have useful results available in these cases,
but they are derived from a deeper path analysis of the L\'evy process and will be presented separately.

As we illustrated in Section \ref{s7}, the formulae listed in Table 1
 can be quite easily calculated in special cases, and it would be very useful to develop computational
and/or simulation methods to deal with more general cases.

\setcounter{equation}{0}\section{Proofs}\label{s9}

\bigskip \noindent {\bf Proof of Theorem \ref{th1}}.\
 Assume \eqref{cramer} and \eqref{mu*}.
Subtracting \eqref{om} from \eqref{tm} and dividing by $P(\tau_u<\infty)$
gives, for $x\ge 0$,
\begin{eqnarray}\label{over1}
P^{(u)}\left(X_{\tau_u}-u>x\right)
&=&
\frac{1}{P(\tau_u < \infty)}
\int_{[0,u]}\pibar_{\cal H}(u+x-y) V_H(\rmd y)
\nonumber\\
&\sim&
q^{-1} \nu_0 m^* e^{\nu_0 u}\int_{[0,u]}\pibar_{\cal H}(u+x-y) V_H(\rmd y)  \quad {\rm (by}\ \eqref{BDmain})
\nonumber\\
&=&
q^{-1} \nu_0 m^* \int_{[0,u]} e^{\nu_0 (u-y)}\pibar_{\cal H}(u+x-y) U^*(\rmd y)\quad {\rm (by}\ \eqref{tauV})
\nonumber\\
&=&
q^{-1} \nu_0 m^* \int_{[0,u]} e^{\nu_0 y} \pibar_{\cal H}(x+y) \rmd_y\left(U^*(u)-U^*(u-y)\right).
\end{eqnarray}
Thus, integrating by parts,
 \begin{eqnarray}\label{over}
P^{(u)}\left(X_{\tau_u}-u>x\right)
&\sim&
q^{-1} \nu_0 m^* e^{\nu_0 u} \pibar_{\cal H}(u+x) U^*(u)-
\nonumber\\
&&
\
q^{-1} \nu_0 m^*
\int_{[0,u]}  \left(U^*(u)-U^*(u-y)\right)
\rmd_y( e^{\nu_0 y} \pibar_{\cal H}(x+y)).
\end{eqnarray}
For later reference, note that for any $x\ge 0$,
 \begin{eqnarray}\label{intfi}
\int_{[0,\infty)} y|\rmd_y\left(e^{\nu_0 y} \pibar_{\cal H}(x+y)\right)|
&\le&
\int_{[0,\infty)} y \nu_0e^{\nu_0 y} \pibar_{\cal H}(x+y)\rmd y+\int_{[0,\infty)} y e^{\nu_0 y} \Pi_{\cal H}(x+\rmd y)\nonumber\\
&=&
\int_{[0,\infty)} \int_{[0,y)} \nu_0 ze^{\nu_0 z} \Pi_{\cal H}(x+\rmd y)\rmd z+\int_{[0,\infty)}y e^{\nu_0 y} \Pi_{\cal H}(x+\rmd y\nonumber)\\
&\le&
2\int_{[0,\infty)}y e^{\nu_0 y} \Pi_{\cal H}(x+\rmd y)<\infty,
 \end{eqnarray}
by \eqref{rid} since $m^*<\infty$.

First assume $\rmd_H>0$.  Then $\rmd_{H^*}=\rmd_H>0$ by \eqref{dpstar}, and so
by Prop III.1 of Bertoin (1996),
$U^*(y)\le c_1 y$ for all $y\ge 0$, for some constant $c_1>0$.
Thus the first term on the RHS of \eqref{over} tends to 0 as $u\to\infty$, by \eqref{picon}. Next by
subadditivity of $U^*$, we have
\begin{equation}\label{subadd}
0\le U^*(u)-U^*(u-y)\le U^*(y)\le c_1 y,\ {\rm for\ all}\
0\le y\le u.
\end{equation}
Thus by \eqref{intfi} we may apply dominated convergence to the second term on the RHS of \eqref{over}.
It then follows from \eqref{U*con} and \eqref{picon} that
\begin{eqnarray*}
\lim_{u\to\infty}P^{(u)}\left(X_{\tau_u}-u>x\right)
&=&
-q^{-1} \nu_0 m^*
\int_{[0,\infty)} \left(\frac{y}{m^*}\right) \rmd_y\left(e^{\nu_0 y} \pibar_{\cal H}(x+y)\right)
\nonumber\\
&=&
q^{-1} \nu_0\int_{[0,\infty)}e^{\nu_0 y} \pibar_{\cal H}(x+y)\rmd y
\nonumber\\
&=&
q^{-1} \int_{[0,\infty)}\left(e^{\nu_0 y}-1\right) \Pi_{\cal H}(x+\rmd y),
\end{eqnarray*}
for all $x\ge 0$, proving \eqref{t2a} when $\rmd_H>0$.

Now assume $\rmd_H=0$.  In that case
we can only deduce
from Prop III.1 of Bertoin (1996) that for each $y_0>0$, there is a constant $c_2=c_2(y_0)>0$ such that
$U^*(y)\le c_2 y$ for $y\ge y_0$, and so
\eqref{subadd} has to be modified accordingly.
To account for this we need to consider separately the cases $x>0$ and $x=0$.  The first term on the RHS of \eqref{over} tends to $0$, as before, in either case. When $x>0$ we can change variable in the second term  to get
\[
q^{-1} \nu_0 m^* e^{-\nu_0 x}
\int_{[x,x+u]} \left(U^*(u)-U^*(u-(y-x))\right)
\rmd_y( e^{\nu_0 y} \pibar_{\cal H}(y)).
\]
By subadditivity, the integrand is bounded above by $U^*(y-x)\le U^*(y)\le c_2y$ for some $c_2=c_2(x)$, if $y\ge x>0$.
Since $y|\rmd_y( e^{\nu_0 y} \pibar_{\cal H}(y))|$
is integrable by \eqref{intfi}, we can
apply dominated convergence to obtain
\begin{eqnarray*}
\lim_{u\to\infty}P^{(u)}\left(X_{\tau_u}-u>x\right)
&=&
-q^{-1} \nu_0 m^* e^{-\nu_0 x}
\int_{[x,\infty)} \left(\frac{y-x}{m^*}\right) \rmd_y\left(e^{\nu_0 y} \pibar_{\cal H}(y)\right)
\nonumber\\
&=&
q^{-1} \int_{[0,\infty)}\left(e^{\nu_0 y}-1\right) \Pi_{\cal H}(x+\rmd y),
\end{eqnarray*}
just as before. It remains to deal with the case $\rmd_H=0$ and $x=0$.
But then \eqref{t2a} follows immediately from \eqref{kes} and \eqref{qid}.

Next we prove \eqref{t5a} and \eqref{t6a} for $x>0$. We first observe that for any $x>0$
\ben\ba
\int_{(x,u]} &e^{\nu_0 y}\pibar_{\cal H}(y)\rmd_y\left(U^*(u)-U^*(u-y)\right)\\
&\qquad=
e^{\nu_0 u} \pibar_{\cal H}(u)U^*(u) - e^{\nu_0 x} \pibar_{\cal H}(x)\left(U^*(u)-U^*((u-x)-)\right)\\
&\qquad\quad  -\int_{(0,\infty)} {\bf 1}_{\{x< y\le u\}}\left(U^*(u)-U^*(u-y)\right)\rmd_y( e^{\nu_0 y} \pibar_{\cal H}(y)).
\ea\een
Now for some constant $c_2=c_2(x)$,
$U^*(u)-U^*(u-y)\le c_2 y$ for $y\ge x>0$, as observed above.
Thus, as $u\to\infty$, the first term tends to $0$ by \eqref{picon},  while
$y|\rmd_y( e^{\nu_0 y} \pibar_{\cal H}(y))|$
is integrable on $(0,\infty)$ by \eqref{intfi}.
So by dominated convergence
\begin{eqnarray}\label{xu}
&&
\int_{(x,u]} e^{\nu_0 y}\pibar_{\cal H}(y)\rmd_y\left(U^*(u)-U^*(u-y)\right)
\nonumber\\
&&\qquad\quad \to - \left(\frac {x}{m^*}\right)e^{\nu_0 x} \pibar_{\cal H}(x)-\int_{(x,\infty)}
\left(\frac {y}{m^*}\right)
\rmd_y( e^{\nu_0 y} \pibar_{\cal H}(y))
\nonumber\\
&&\qquad\quad =\frac {1}{m^*}\int_{(x,\infty)}e^{\nu_0 y}\pibar_{\cal H}(y)\rmd y.
\end{eqnarray}

Now consider \eqref{t6a} for $x>0$.
Subtract \eqref{3.7} from \eqref{tm}  to obtain, for $0\le x\le u$,
\begin{eqnarray*}
P^{(u)}\left(u-\overline X_{\tau_u-}>x\right)
 &=&
\frac{1}{P\left(\tau_u<\infty \right)}
\int_{[0,u-x)} \overline{\Pi}_{\cal H}(u-y) V_H(\rmd y)\nonumber\\
&\sim&
q^{-1} \nu_0 m^* e^{\nu_0u}  \int_{[0,u-x)}
\pibar_{\cal H}(u-y) e^{-\nu_0y}U^*(\rmd y)
\nonumber\\
 &=&
  q^{-1} \nu_0 m^* \int_{(x,u]} e^{\nu_0 y}\pibar_{\cal H}(y)\rmd_y\left(U^*(u)-U^*(u-y)\right)
 \nonumber\\
 &\to&
\frac {\nu_0}{q}\int_{(x,\infty)}e^{\nu_0 y}\pibar_{\cal H}(y)\rmd y
\nonumber\\
\end{eqnarray*}
by \eqref{xu}. Since
\ben
\frac {\nu_0}{q}\int_{(0,\infty)}e^{\nu_0 y}\pibar_{\cal H}(y)\rmd y=
\frac{1}{q}\int_{(0,\infty)}\left(e^{\nu_0 y}-1\right) \Pi_{\cal H}(\rmd y)=1-\frac {\nu_0\rmd_H}{q}
\een
by \eqref{qid}, \eqref{t6a} then follows.

To prove \eqref{t5a} for $x>0$,
subtract \eqref{um} from \eqref{tm} to obtain
\begin{eqnarray} \label{A}
P\left(u-X_{\tau_u-} > x, \tau_u<\infty\right)
&=&
\int_{[0,u]}\overline{\Pi}_{\cal H}(u-y) V_H(\rmd y)
-\int_{y\in[u-x, u]} g_x(u-y)V_H(\rmd y).
 \nonumber\\
\end{eqnarray}
The first term on the RHS of  \eqref{A},
divided by $P(\tau_u<\infty)$, is just $P^{(u)}\left(X_{\tau_u}-u > 0\right)$, by the first equality in \eqref{over1}. Applying  the already proven \eqref{t2a} with $x=0$, we get
 \begin{eqnarray}\label{term1}
\lim_{u \to \infty}\int_{[0,u]}\overline{\Pi}_{\cal H}(u-y) V_H(\rmd y)=\lim_{u \to \infty}P^{(u)}\left(X_{\tau_u}-u > 0\right)
=1-\frac{\nu_0\rmd_H}{q}.
\end{eqnarray}
The second term in \eqref{A}, when divided by $P(\tau_u<\infty)$, is
\begin{eqnarray}\label{9y}
 \frac{1}{P(\tau_u<\infty)}\int_{[u-x, u]} g_x(u-y)e^{-\nu_0y}U^*(\rmd y)
&\sim&
\frac{\nu_0m^*}{q}\int_{[0, x]} g_x(y)e^{\nu_0y}\rmd_y\left(U^*(u)- U^*(u-y)\right). \nonumber\\
\end{eqnarray}
Take $\delta\in (0,x)$ and write $\int_{[0,x]}=\int_{[0,\delta]}+ \int_{(\delta,x]}$.
For the integral over $[0,\delta]$ we have
\be\ba\label{smd}
\int_{[0, \delta]} g_x(y)e^{\nu_0y}\rmd_y\left(U^*(u)- U^*(u-y)\right)
&\le \int_{[0, \delta]} \pibar_{\cal H}(y)e^{\nu_0y}\rmd_y\left(U^*(u)- U^*(u-y)\right)\\
\ea\ee by \eqref{ginf}.
On the RHS of \eqref{smd}, write
$\int_{[0,\delta]}= \int_{[0,u]}-\int_{(\delta,u]}$, then
use \eqref{over1} and \eqref{t2a}, both  with $x=0$, and  \eqref{xu} with $x=\delta$, to see that,
as $u\to\infty$, the RHS of  \eqref{smd} converges to
\be\ba\label{dellim}
\frac{1}{\nu_0m^*} \int_{[0,\infty)}\left(e^{\nu_0 y}-1\right) \Pi_{\cal H}(\rmd y)- \frac{1}{m^*}\int_{(\delta,\infty)}e^{\nu_0 y}\pibar_{\cal H}(y)\rmd y
&=\frac{1}{m^*} \int_{[0,\delta]}e^{\nu_0 y}\pibar_{\cal H}(y)\rmd y.
\ea\ee
For the integral over $(\delta, x]$ in \eqref{9y}, integration by parts gives
\be\ba\label{2t}
\int_{(\delta, x]} g_x(y)&e^{\nu_0y}\rmd_y\left(U^*(u)- U^*(u-y)\right)\\
&=g_x(x)e^{\nu_0x}\left(U^*(u)- U^*((u-x)-)\right) - g_x(\delta)e^{\nu_0\delta}\left(U^*(u)- U^*((u-\delta)-)\right)\\
&\qquad\qquad-\int_{(\delta, x]}\left(U^*(u)- U^*(u-y)\right)\rmd_y\left(g_x(y)e^{\nu_0y}\right).
\ea\ee
We can apply bounded convergence to this integral since $|\rmd_y\left(g_x(y)e^{\nu_0y}\right)|$
is a finite measure on $(\delta, x]$, and the integrand is bounded because
$U^*(u)- U^*(u-y)\le U^*(y)\le U^*(x)$
for $0\le y\le x$.  Thus
\be\ba\label{2tt}
\int_{(\delta, x]} g_x(y)&e^{\nu_0y}\rmd_y\left(U^*(u)- U^*(u-y)\right)\\
&\to \frac{x}{m^*}g_x(x)e^{\nu_0x} - \frac{\delta}{m^*}g_x(\delta)e^{\nu_0\delta}
-\int_{(\delta, x]}\left(\frac y{m^*}\right) \rmd_y\left(g_x(y)e^{\nu_0y}\right)\\
&=\frac {1}{m^*}\int_{(\delta, x]}g_x(y)e^{\nu_0y} \rmd y,
\ea\ee
after another integration by parts.
Combining \eqref{A}--\eqref{dellim}  and \eqref{2tt}, and letting $\delta\to 0$, proves \eqref{t5a} for $x>0$.

To complete the proof we must show \eqref{t5a} and \eqref{t6a} hold for $x=0$,
that is, the last two equalities of \eqref{t3a} hold.
By (5.2) of Griffin and Maller (2011), for any L\'evy process,
$
P(X_{\tau_u-}<u=X_{\tau_u}, \tau_u<\infty)=0.
$
Since
$
X_{\tau_u-}\le u
$
on $\{\tau_u<\infty\}$, it follows that
$
P(X_{\tau_u}=u, \tau_u<\infty)\le P(X_{\tau_u-}=u, \tau_u<\infty).
$
On the other hand
$
X_{\tau_u-}\le \Xbar_{\tau_u-}\le u
$
on $\{\tau_u<\infty\}$, thus we may conclude that
\ben
P^{(u)}\left(X_{\tau_u}=u\right)\le
P^{(u)}\left(X_{\tau_u-}=u\right)\le P^{(u)}\left(\overline X_{\tau_u-}=u\right).
\een
Since \eqref{t2a} has already been proved for $x=0$, we have
\be
\lim_{u\to\infty} P^{(u)}\left(X_{\tau_u}=u\right)
\ge
 \frac{\nu_0\rmd_H}{q},
 \ee
while for every $x>0$,
\be\ba
\limsup_{u\to\infty} P^{(u)}\left(\overline X_{\tau_u-}=u\right)
&\le \lim_{u\to\infty} P^{(u)}\left(u- \overline X_{\tau_u-}\le x\right)\\
&= \frac{\nu_0\rmd_H}{q}+ \frac{1}{q}\int_0^x \pibar_{\cal H}(y)\rmd \left(e^{\nu_0 y}\right),
\ea\ee
by \eqref{t6a}.
Letting $x\to 0$ completes the proof.\hfill\halmos


\begin{thebibliography}{99}

\bibitem{a02}
Asmussen, S. (2002) {\em Applied Probability and Queues}, 2nd
Edition.  Wiley, Chichester.

\bibitem{Bert}
Bertoin, J. (1996). {\em L\'evy Processes.} Cambridge Univ.  Press.

\bibitem{BD94}
Bertoin, J. and Doney, R.A. (1994)
Cram\'er's estimate for L\'evy processes.
Statistics and  Probability Letters, {\bf 21}, 363--365.

\bibitem{BD96}  Bertoin, J. and Doney, R.A. (1996)
Some asymptotic results for transient random walks.
Advances in Applied Probability,  \textbf{28} 207--226.

\bibitem{BK}
Biffis, A. and Kyprianou, A.E. (2010)
A note on scale functions and the time value of ruin for L\'evy insurance risk processes, Insurance: Mathematics and Economics 46, 85-91.


\bibitem{BGT}
Bingham, N.H., Goldie, C.M. and Teugels, J.L. (1987). {\em
Regular Variation.} Cambridge University Press, Cambridge.

\bibitem{CNW}
Chover, J., Ney, P. and Wainger, S. (1973). Degeneracy properties of subcritical branching processes. Annals of Probability, {\bf 1}, 663--673.

\bibitem{DK}
Doney, R.A. and Kyprianou, A. (2006). Overshoots and undershoots
of L\'evy processes. Annals of  Applied Probability, {\bf 16}, 91--106.

\bibitem{DM}
Doney, R.A. and Maller, R.A. (2004). Moments of passage times for L\'evy processes.
Annals of the Institute Henri Poincar\'e, Probability and  Statistics, {\bf 40}, 279--297.

\bibitem{EG}
Embrechts, P. and Goldie, C.M. (1982). On convolution tails.
Stochastic  Processes and  Applications, {\bf 13}, 263--278.

\bibitem{EV}
Embrechts, P. and Veraverbeke N. (1982).
Estimates for the probability of ruin
with special emphasis on the possibility
of large claims, Insurance: Mathematics and Economics I, 55-72.


\bibitem{EGV}
Embrechts, P., Goldie, C.M. and Veraverbeke N. (1979).
Subexponentiality and infinite divisibility. Zeitschrifte fur
Wahrscheinlichkeitstheory Verwe Gebiete, {\bf 49} 335--347.

\bibitem{grandell}
Grandell, J. (1991). {\em Aspects of Risk Theory.} New York: Springer-Verlag.

\bibitem{GM1} Griffin, P.S. and Maller, R.A. (2010)
Path decomposition of ruinous behaviour for a general L\'evy insurance risk process.  To Appear in Annals of Applied Probability.

\bibitem{GM2}
Griffin, P.S. and Maller, R.A.  (2011) The time at which a L\'evy process creeps, submitted.

\bibitem{hk}
Hubalek, F. and Kyprianou, A. E. (2010). Old and new examples of scale
functions for spectrally negative L\'evy processes. To appear in: 6th Seminar on
Stochastic Analysis, Random Fields and Applications, Eds R. Dalang, M. Dozzi, F. Russo. Progress in Probability, Birkhauser.

\bibitem{kl}
Kl\"{u}ppelberg, C. (1989). Subexponential distributions and characterizations of related classes.
Probability Theory and Related Fields, {\bf 82}, 259--269.

\bibitem{kkm}
Kl\"{u}ppelberg, C., Kyprianou A. and Maller, R.A. (2004).
Ruin probability and overshoots for general L\'{e}vy insurance risk processes.
Annals of  Applied Probability, {\bf 14}, 1766-1801.

\bibitem{kyp}
Kyprianou A. (2005).
{\em Introductory Lectures on Fluctuations of L\'{e}vy
Processes with Applications}. Springer: Berlin Heidelberg New York.

\bibitem{PM}
Park, H.S. and Maller, R.A. (2008)
Moment and MGF convergence of  overshoots and undershoots
for L\'{e}vy insurance risk processes, Advances in Applied Probability, 40, 716--733.

\bibitem{sato}
Sato, K. (1999). {\em L\'evy Processes and Infinitely Divisible Distributions.}
Cambridge University Press, Cambridge.


\bibitem{schmid}
Schmidli, H. (1995) Cram\'er-Lundberg approximations for ruin probabilities of risk processes perturbed by diffusion, Insurance: Mathematics and Economics,  16, 135-149.

\bibitem{WangLei}
Tang, Q. and Wei, L. (2010)
Asymptotic aspects of the Gerber-Shiu function in the renewal risk model using Wiener-Hopf factorization and convolution equivalence, Insurance: Mathematics and Economics 46, 19-31.

\bibitem{Vig}
Vigon, V. (2002). Votre L\'{e}vy rampe-t-il?
Journal of the  London Mathematics Society, {\bf 65}, 243--256.

\end{thebibliography}
\end{document}